\documentclass[12pt]{article}

\usepackage{hyperref}
\usepackage{graphicx}
\usepackage[ansinew]{inputenc}
\usepackage{enumerate}
\usepackage{amsmath}
\usepackage{amsfonts}
\usepackage{amssymb}
\usepackage{amsthm}
\usepackage{multirow}
\usepackage{graphicx}
\usepackage{color}
\usepackage{microtype}
\input xy
\xyoption{all}
\usepackage{url}
\usepackage{varioref}
\usepackage{eso-pic}
\usepackage{tikz}
\usepackage{float}
\usepackage{dsfont}
\usepackage{listings}

\lstdefinelanguage{Magma}
{
keywords={for,end,if,then,else,elif,while,function,return,cat,&,and,or },
morekeywords={Seqset,Setseq,Polytope,AutomorphismGroup,RowSequence,IdentifyGroup,
	      Subgroups,PermutationMatrix,Generators,MatrixGroup,Transpose},
sensitive=false,
morecomment=[l]{//},
morecomment=[s]{/*}{*/},
morestring=[b]",
}

\lstnewenvironment{codice_magma}[1][]
{\lstset{basicstyle=\small\ttfamily, columns=fullflexible,
language=Magma,
keywordstyle=\color{red}\bfseries,
commentstyle=\color{blue},
numbers=left, numberstyle=\tiny,
stepnumber=2, numbersep=5pt, frame=single, #1}}{}

\usepackage{multicol}
\theoremstyle{definition}

\theoremstyle{remark}

\begin{document}

\title{\textbf {\textit{MadeinMath}. \\ A Math Exhibition out of Many Threads}}
\author{Gilberto Bini \thanks{Universit\`a degli Studi di Milano, Dipartimento di Matematica "Federigo Enriques" and {\it matematita} - Centro interuniversitario di ricerca per la comunicazione e l'apprendimento informale della matematica (\url {www.matematita.it}),
Via Cesare Saldini 50, 20133 Milano, Italy, E-mail: {\tt gilberto.bini@unimi.it}}}

\date{\today}

\maketitle



\begin{abstract}
We describe some motivations lurking behind the making of a Math exhibition. In particular, we refer to \textit{MadeInMath} and its two set-ups in Triennale di Milano and Museo della Scienza in Trento (MUSE), where the life and works of Italian mathematicians played an important role, like Guido Castelnuovo, Federigo Enriques and Vito Volterra.
\end{abstract}

\section{Introduction}\label{intro} 

Very often, those who try to solve Math problems have the impression of being in a maze: students struggling with their studies, architects intent to determine geometric configurations, researchers of mathematics making an attempt to solve - among many others - the problem of the millennium. It is easy to get lost in a maze because you do not have access to a map that shows you the way, or to a magic wand that allows you to see from a different perspective the maze of streets in front of you and lead to the way out. As in the famous myth of the labyrinth and the Minotaur, it may come handy to you a \textit{wire},  the legendary Ariadne's wire, which allows you to find the high road if the wire is slowly disentangled.

Undoubtedly, the maze is persistent and indestructible in human imagination! And so is the myth of Ariadne and Theseus, where hope, love and desire for salvation are interwoven. Very often, the myths concerning mazes (including that of Palace Cnosso, which is certainly the most famous) are connected to the symbol of the {\it spiral} (a kind of toy-maze). As told in some fairy tales \cite{P}, in the symbolism of the spiral a preteen is lost, overcomes difficult trials, has magical aid, finds enemies, etc., but eventually learns and comes back as an adult. In a sense, he dies and is reborn. 

Mazes and wires, fairy tales and spirals may look completely disjoint from a mathematical context. To the contrary, the metaphor of the labyrinth and the sense of loss seem quite natural when faced with research problems - at various levels. And as a wire in a labyrinth can be useful not to remain in the grip of anguish of not knowing how to deal with the situation, so a metaphorical version of it is very useful for researchers, whether they are students dealing with homework assignments or professional mathematicians dealing with research problems. Indeed, a guidance helps finding the possible solutions via models based on a profound mix of theory and practice. Equally important is the metaphor between spiral and learning because both do not come out of a linear progression of a succession of knowledge, but rather a circular development of skills that are perfected over time thanks to the constant references between one phase and the other of the growth.

The starting point in these cases lies in the already formalized and abstract language
of mathematics. Abstraction in mathematics is not an escape from reality, but an attempt to understand it better: sometimes it is useful to move away a little, like when we view things at the right distance. From there, reality appears to us in all its essential features, and situations that seem to have little or nothing to do with each other can turn out to be similar, or even identical. When, beyond different details and appearances, their behaviour is determined by the same rules, mathematics say that they have the same \textit{structure}. Then it makes sense to study these
structures in themselves, forgetting the more or less concrete situations that might have generated them. Abstraction is a tool that various areas of research have in common, both those that begin with concrete problems and those that originate within mathematics. In a sense, abstraction makes studying the circulation of blood similar to studying the dynamics of the machines that fill Tetra Pak packages; indeed, both have been seen to obey the same rules. Abstraction is classifying the different algebraic surfaces on the basis of certain parameters that are considered relevant.

How can we put into practice these metacognitive observations? How can we raise our interlocutors from the sense of loss that many feel when dealing with mathematics? And give them some gentle nudges (metaphorical wires) to disentangle the key concepts of Mathematics from the spiraling unfolding of real phenomena? And ultimately recover the platonic nature of Math?

Possible answers to this set of questions are not so easy to be formulated, but an attempt in this direction was made within the math exhibition \textit{MadeInMath}. The first staging of this exhibition was developed at the Triennale di Milano from 14 September 2014 to 23 November 2014 with more than 40,000 visitors. Curated by Gilberto Bini, Maria Ded\`o, Simonetta Di Sieno (Center "matematita", University of Milan) and Renato Betti, Angelo Guerraggio (Centro PRISTEM, University L. Bocconi), the exhibition was organized by the Center "matematita" and the Centro PRISTEM in collaboration with Dr. Vincenzo Napolano (INFN) and sponsored, among others, by many Mathematics Departments of the major Italian universities. The second set up of the exhibition was made at {\it Muse} - Museum of Science (Trento) from 26 February to 26 June 2016. Marco Andreatta (CIRM-FBK, Muse and Department of Mathematics, Trento) was an additional curator of this second set-up.

\section{Out of Many Threads}\label{sec1}

As suggested by the title of this article, the format of the exhibition includes the presentation of various mathematical themes all related to one another, such as mathematical modeling, games and puzzles, the concept of numbers, the visualization of the fourth dimension, etc. These various topics may seem at first disconnected from each other, and some in fact are. The entrance to the exhibition looks a little like the entrance in a maze. In order not to get caught by the sense of disorientation and to appreciate the exhibition, the visitor is invited through various hints, such as multimedia installations, written texts and animations, etc., to choose a possible `wire` in order to tackle the course and see - occasionally - connections with different themes, as long as he/she has time and desire. 

{\it So the exhibition is more modeled on a spiral-like format rather than a linear-like format in the sense that the visit is spread along a path which does not follow necessarily the units of space and time.} 

As an example, the visitor may start with models used for the reproduction of architectural works (like the Duomo in Milano); move on to models used in sports to improve the performance of athletes and those to regulate traffic flows (see Figure 1), such as the study of fluid dynamics for the creation of swimming costumes to compete at the highest levels; go back in time and appreciate the works of the Italian mathematician Vito Volterra (1860-1940), who has opened the way for the application of mathematics to biological sciences with the prey-predator model. 

\begin{figure}[ht]
\centering
\includegraphics[width=10cm, angle=-270]{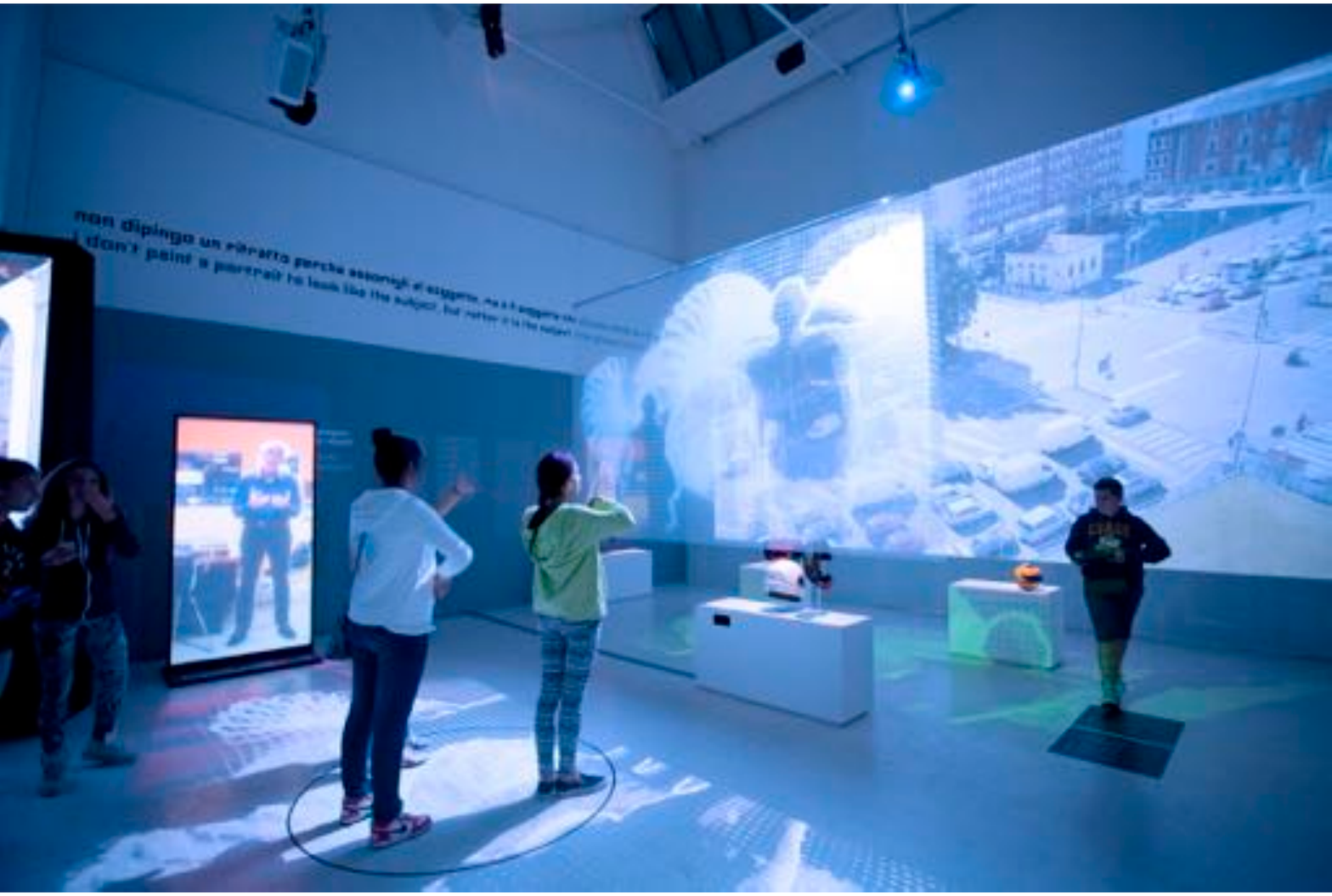}
\caption{{\bf Models from real life (photo courtesy of Carla Mondino).}}
\end{figure}

This `wire` is based on the leitmotif of mathematical models but touches on history of mathematics, applied mathematics, projective geometry and so on. Alternatively, one can select various contributions of some Italian mathematicians, e.g., those from the beginning of the XXth century, and tell the mathematical notions developed through their life: from Guido Castelnuovo (1865-1952) to Federigo Enriques (1871-1946) from Bruno de Finetti (1906-1985) to Tullio Levi-Civita (1873-1941), etc. {\it MadeInMath} is also an exhibition aimed at the present and the future, rich in historical references that serve to better understand how Mathematics is evolving and will evolve in our times. Two installations are dedicated to the present and the future. Finally, some young researchers, who talk about their dreams and their hopes, give meaning to the concerns expressed by the {\it Mappa dei matematici italiani all'estero} (Maps of the Italian Mathematicians abroad) designed by the Unione Matematica Italiana \cite{map} and give substance to the description of the places where research is done and taught.

\begin{figure}[ht]
\centering
\includegraphics[width=10cm, angle=270]{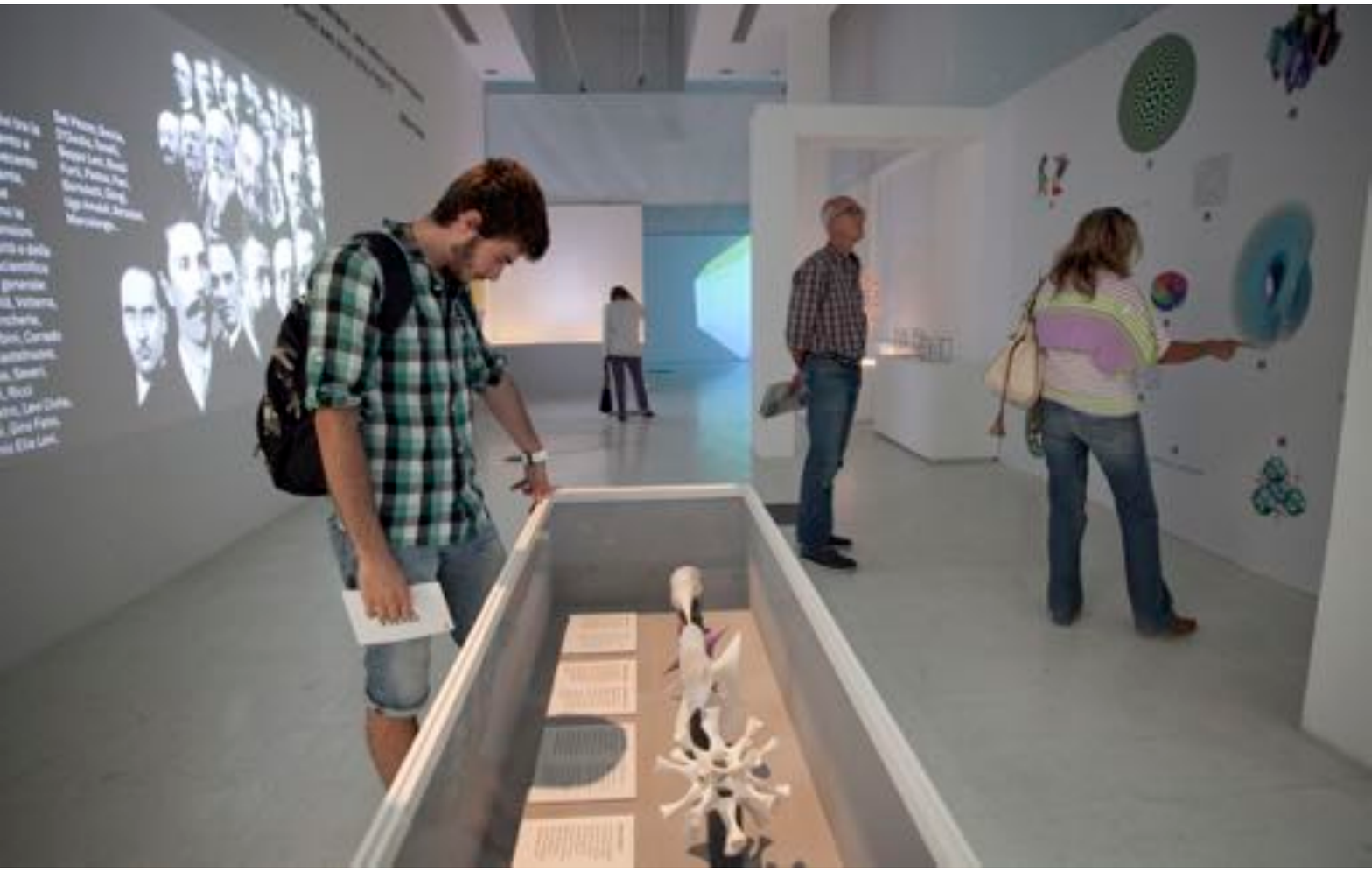}
\caption{{\bf Entering the labyrinth of the exhibition (photo courtesy of Carla Mondino).}}
\end{figure}

\section{An Installation: Algebraic Surfaces} \label{installation} One installation, among others, which has had a great success, has touched on the issue of algebraic surfaces and their relevance in the history of Italian Mathematics. It is inspired to the following excerpt taken from \cite{GC}:

\begin{center}
\textit{
"We had created" in an abstract sense, of course" a large number of models of surfaces in our space or in higher spaces; and we had split these models, so to speak, between two display windows. One contained regular surfaces for which everything proceeded as it would in the best of all possible worlds; analogy allowed the most salient properties of plane curves to be transferred to these. When, however, we tried to check these properties on the surfaces in the other display, that is on the irregular ones, our troubles began, and exceptions of all kinds would crop up ... With the aforementioned procedure, which can be likened to the type used in experimental sciences, we managed to establish some distinctive characters between the two surface families ..."
}
\end{center}

In this passage, abstract objects are turned into models; classifying them is described as an action of disposing them on different `shelves` according to their properties. 

\begin{figure}[ht]
\centering
\includegraphics[width=10cm, angle=-90]{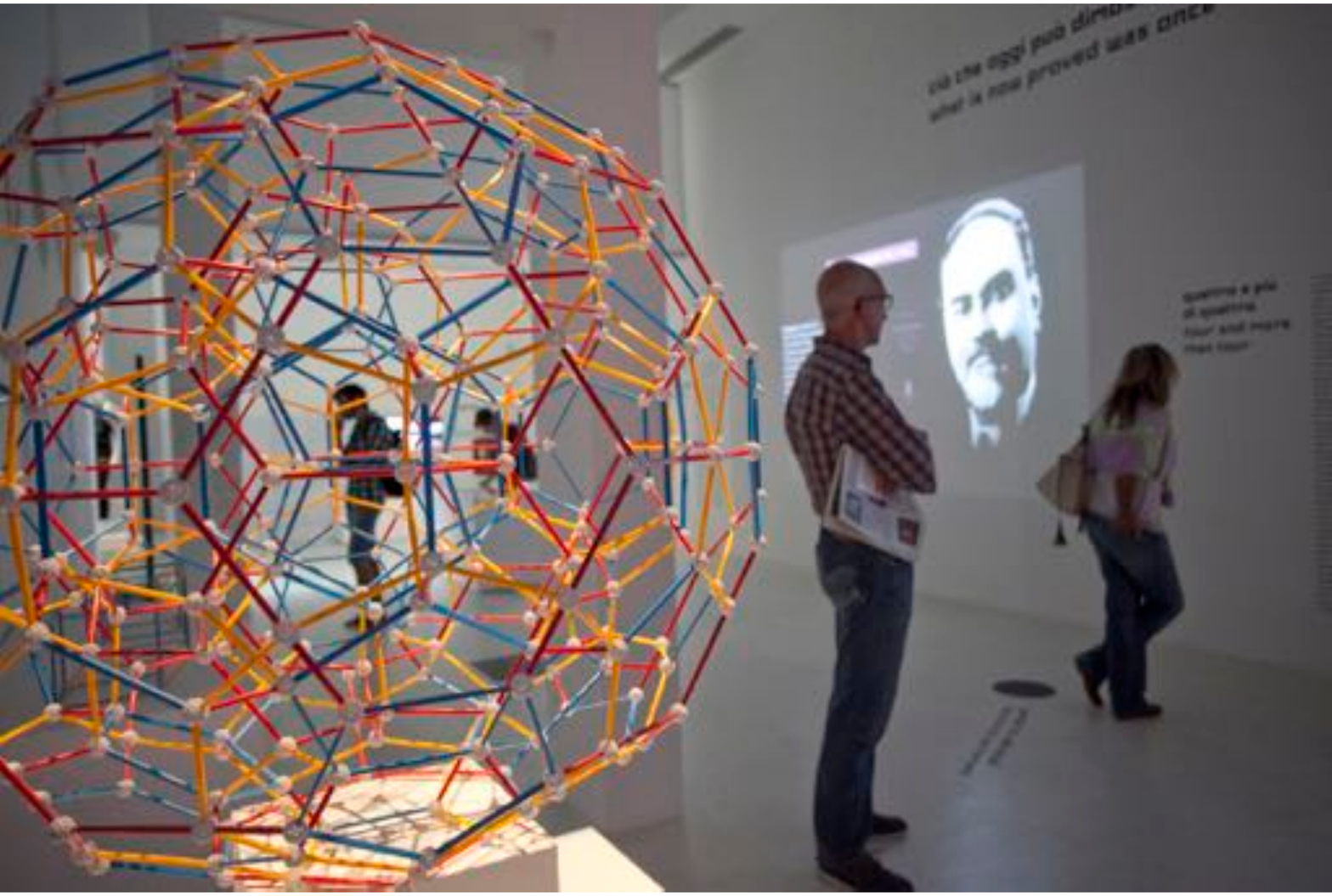}
\caption{{\bf Federigo Enriques from a different perspective (photo courtesy of Carla Mondino).}}
\end{figure}

\pagebreak 

The rendering of this intellectual effort finalized three exhibits. The first was a large wall projection area (see Figure 2 and Figure 3) with some videos concerning famous Italian mathematicians from the 20th century, such as Castelnuovo, Enriques, de Finetti, etc. Opposite to the wall there was another installation with some images and formulas, such as the definition of the homology groups, the composition table of a group of order six, the Fourier transform, Euler's variational equation, etc. These images were briefly explained in order to give an idea of the kind of tools used in doing research. It all seemed weird to the average visitor but the effort was really appreciated because the feedback was positive. People were rather inquisitive. They asked a lot about the intriguing life of the famous mathematicians as well as about the formulas they were presented. The majority of them asked for a down-to-earth explanation of  the most abstract concepts. As explained in the next section, the animators did really a good job in keeping high the attention of visitors without distorting the mathematics that was being presented to them. Finally, close to these two installations, there was a shelf with five 3D models.

\subsection{Some Examples of Surfaces}\label{surfaces}

\begin{figure}[ht]
\centering
\includegraphics[width=5.7cm,angle=270]{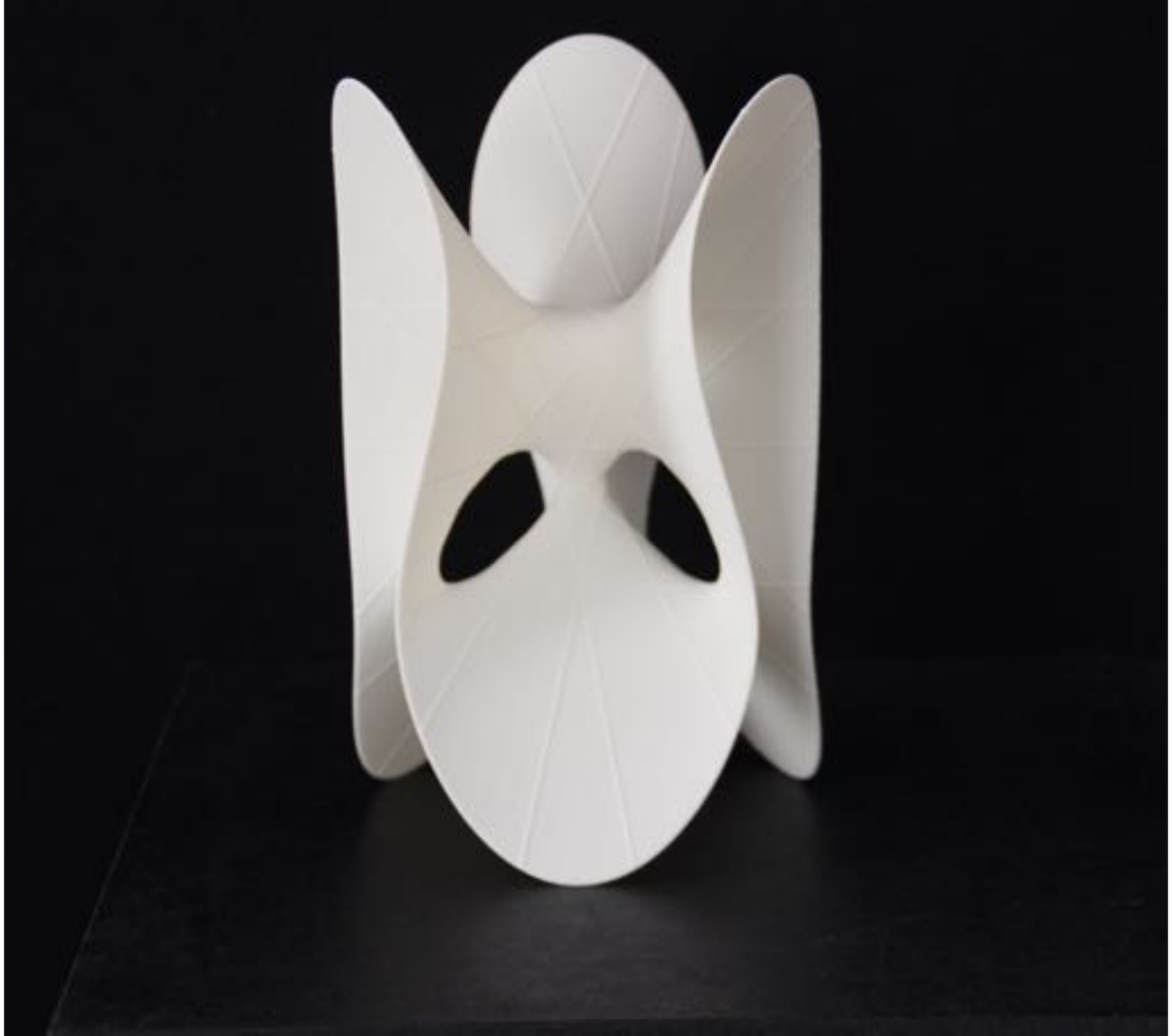} \hskip -0.7cm \includegraphics[width=5.6cm, height=7cm, angle=270]{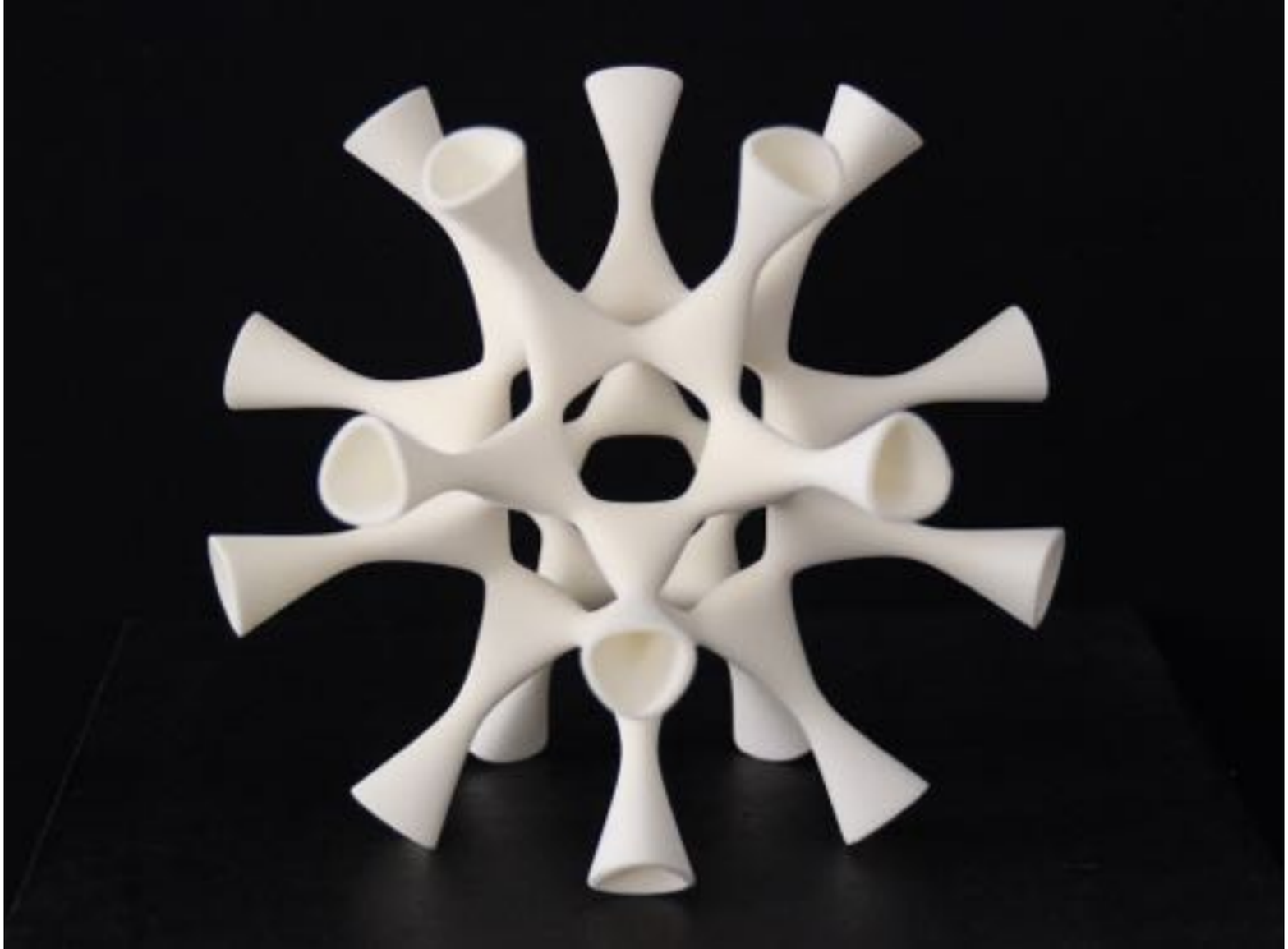} \vskip 0.1cm \hskip 0.65cm
\includegraphics[height=6.9cm,angle=270]{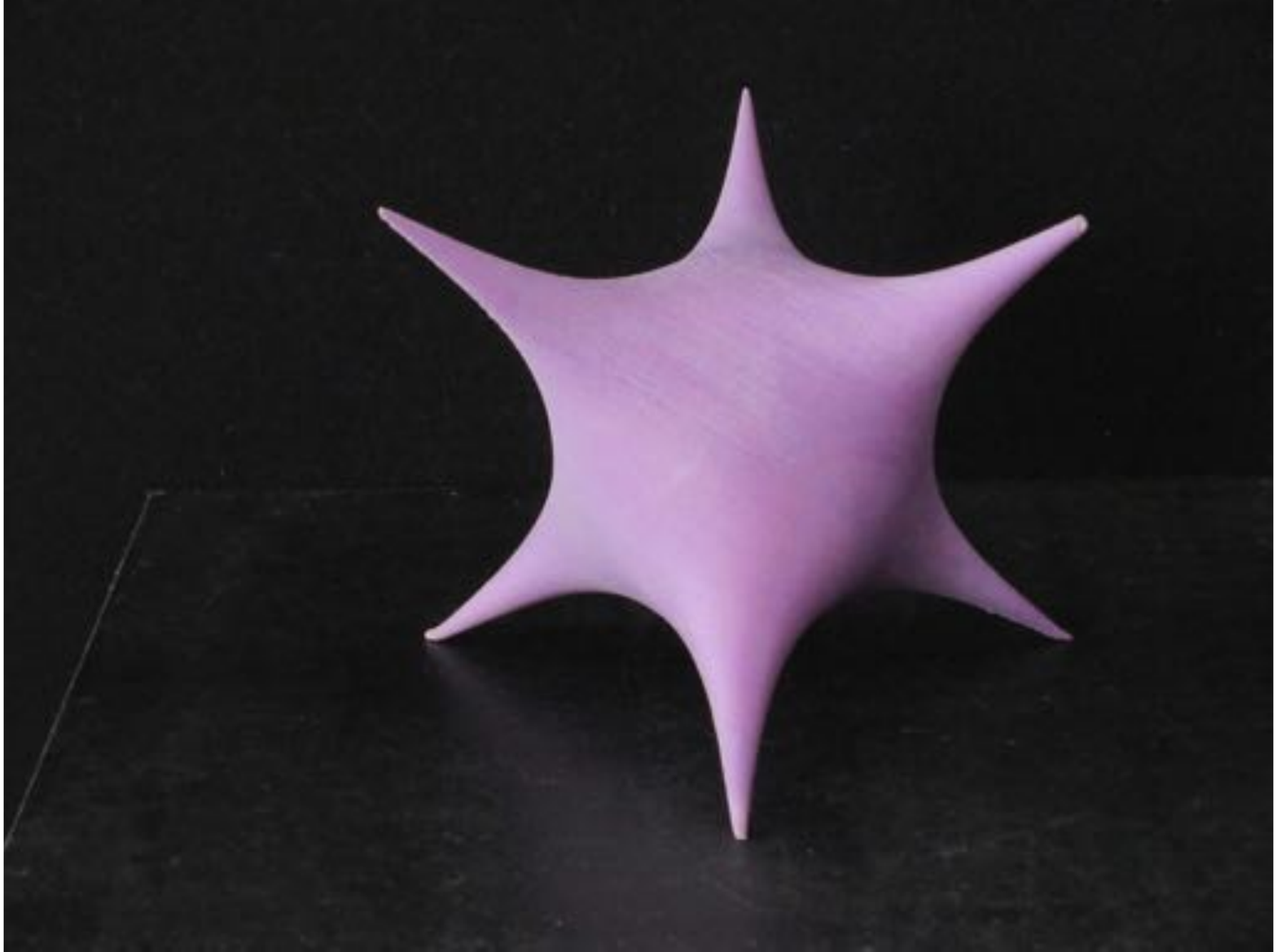}  
\hskip -0.1cm
\includegraphics[width=5cm, height=7cm, angle=270]{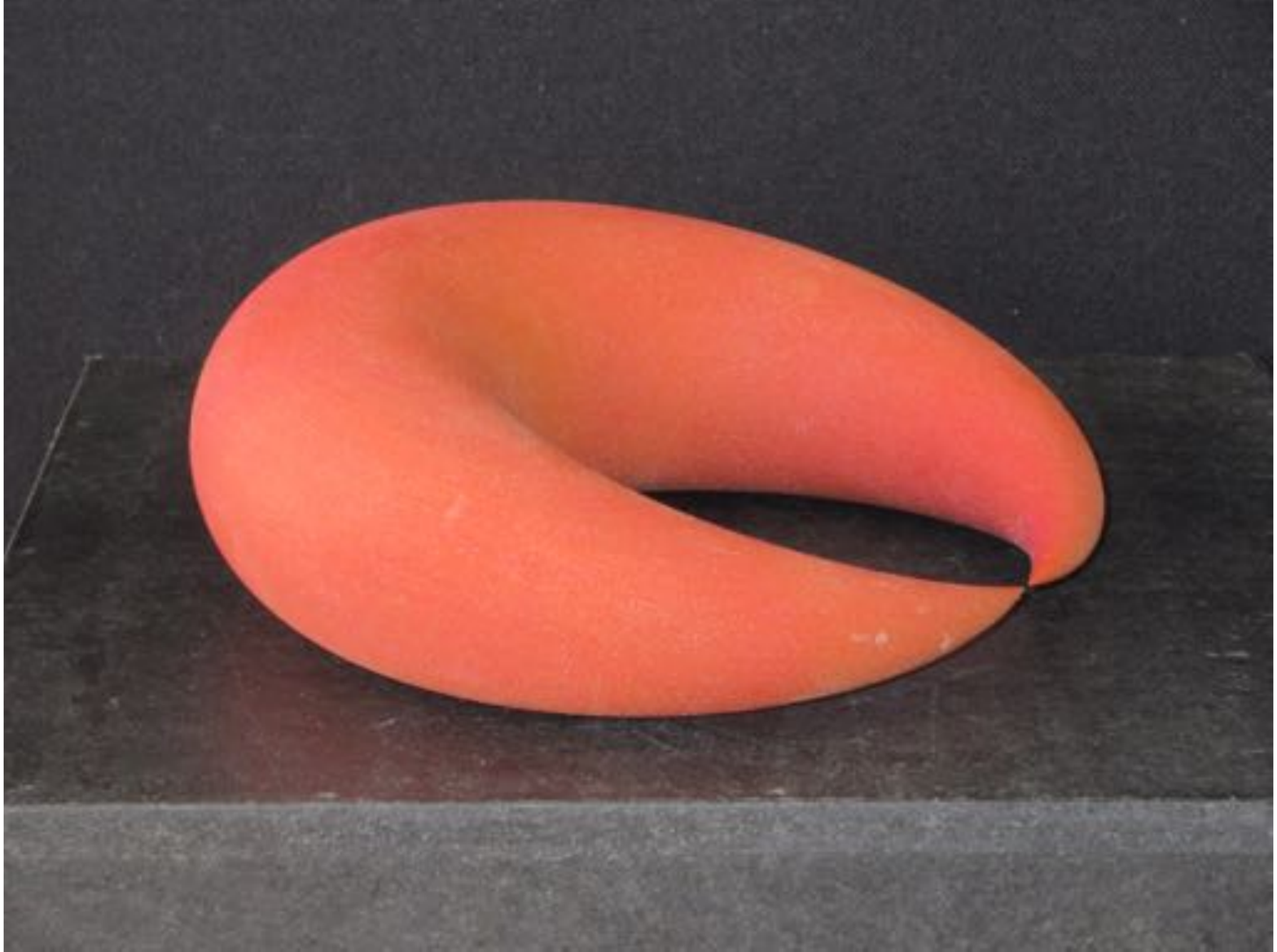}
\caption{{\bf The Clebsch surface and the lines on it (top left); the Barth sextic (top right); a spiky surface (bottom left); a pinched torus (bottom right) - photo courtesy of Enrico Barbanotti.}}
\label{figure}
\end{figure}

\begin{figure}[ht]
\centering
\includegraphics[width=7cm, angle=270]{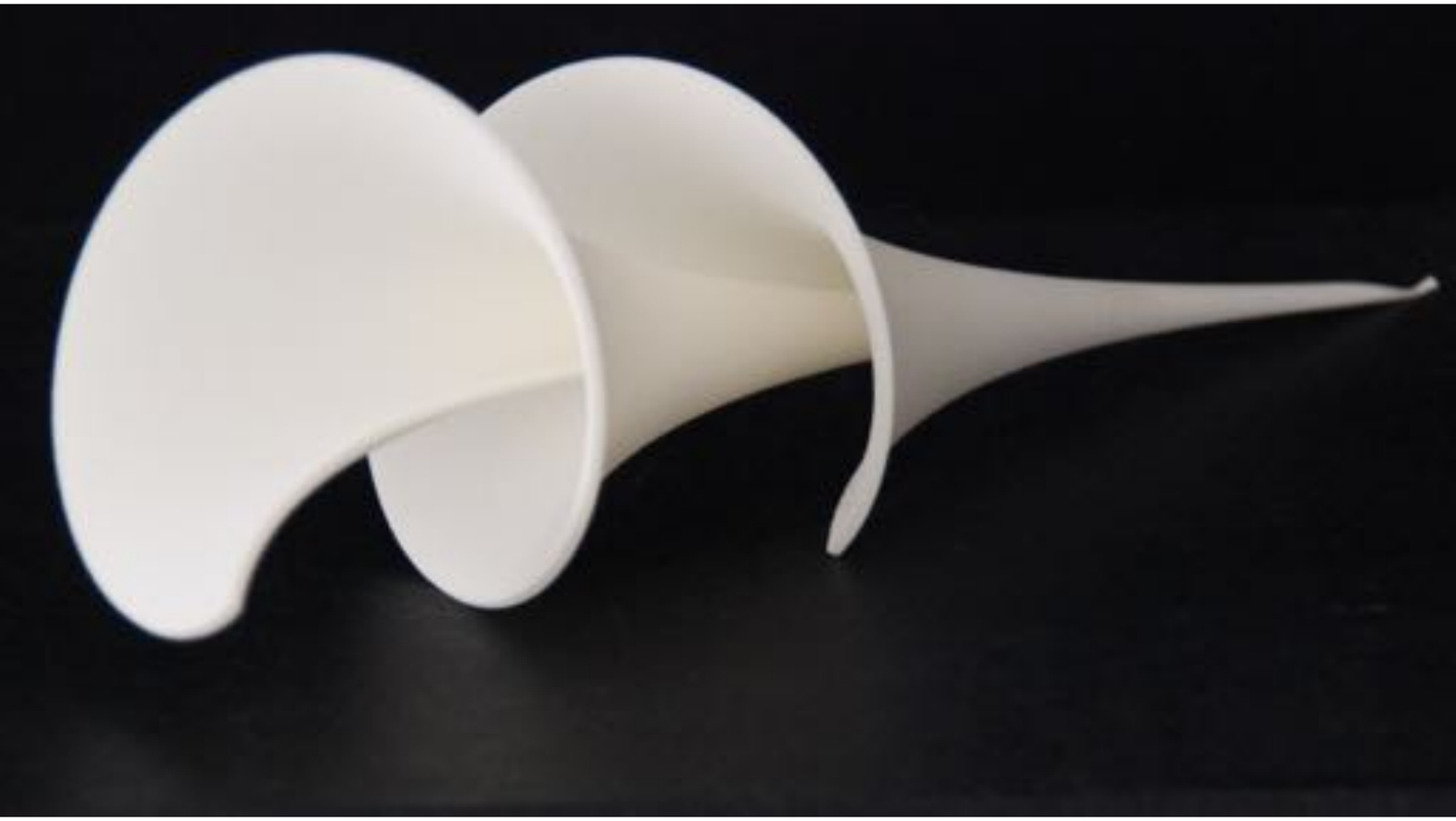}
\caption{{\bf The spiraling Dini surface (photo courtesy of Enrico Barbanotti).}}
\label{dini}
\end{figure}

The Clebsch surface or Klein's icosahedral cubic surface is a non-singular surface  studied by Clebsch (1871) and Klein (1873). This surface contains $27$ lines that can be defined over the real numbers. It is possible to choose the degree three polynomials so that the lines can be seen in three dimensional real space. In the 3D model shown at the exhibition, you could indeed see them (see Figure \ref{figure}). Another surface was one of the Barth surfaces discovered by the German mathematician Wolf Barth in 1996. The chosen one is defined by a degree 6 polynomial in two variables over the complex numbers. The resulting surface is of the complex nodal has 65 ordinary double points: see Figure \ref{figure}. Another surface has six spikes: we can obtain it as a section of an algebraic surface by a polynomial of degree six. The six spikes are situated in the centres of the faces of a cube (see Figure \ref{figure}). There are also surfaces with 4, 8, 12 or 20 spikes that correspond to the centers of the faces of a tetrahedron, octahedron, dodecahedron and icosahedron. Finally, you could find on the shelf the Dini surface and a pinched torus. The former (see Figure \ref{dini}) was discovered by Ulisse Dini (1845-1918) and gives a model of a non-Euclidean geometry. It is far from being algebraic. The latter (see Figure \ref{figure}) is obtained from a doughnut by pinching it. It is not even a topological surface! 

\pagebreak

For further images of algebraic surfaces, see, for instance, \cite{mate}, \cite{imaginary} or \cite{todfer}. All the mathematical notions behind these geometric objects did not in fact reach the average visitor. Nonetheless, they served as a good starting point in order to talk about the interconnection of mathematics with other disciplines, or to talk about a prominent figure like that of Guido Castelnuovo, his life (see, for instance, \cite{boc} for a short bio) and his works \cite{lincei}.

\section{The visitors' impressions} \label{sec4} Multimedia installations and interwoven threads were not the only factors of success for the first set-up of the exhibition \textit{MadeInMath}. The animators were of crucial importance because they were mediators between organizers and visitors. Master students in various degrees in science (Mathematics, Biology, Physics, etc.) have been trained on the topics covered by the exhibition to accompany the many school-groups, but also to follow individual visitors, or small groups of hesitant visitors. For animators it was undoubtedly a formative moment for their training and their personal growth. Each of them, in fact, was free to propose their reading of the contents of the exhibition because the exhibition has various interpretations and visits had different cuts. People reported that they visited it twice and had different impressions each time. The ability to create connections between seemingly unrelated issues is one of the most powerful activity of mathematics, which gave the show an edge. It was thanks to all these facets that each visitor could grasp interesting aspects depending on their age, their training and their own interests. As commented in \cite{tg} by one of the animators, Maurizio Giaffredo, each of us animators took away a wealth of meetings and reflections whose value derives from the diversity of people and ideas that we had the chance to meet. In just two weeks the exhibition was visited by over 40,000 visitors. Of these, about 18,000 were students, mostly from middle and  high school. This means that more than half of the visitors fall into the category `general public`, which was quite a success. Beyond the figures, it was impressive to experience a high level of satisfaction by visitors. An initial indicator in this regard was to see students on a visit with the class during the week and then see them again at the weekend with their parents or other friends to illustrate some of the content of the exhibition as if they were young guides. A second indicator comes from the evidence left by visitors. Let us quote some of the thousands that people delighted in leaving. One of them says: "We knew the beauty of mathematics, but we did not suspect that  it was so beautiful. Thank you also for the love and attention to the teacher's figure". In fact, part of the exhibition gave space to the admirable work done every day by teachers, in this case those of mathematics. We were particularly pleased that this aspect has been apprehended and appreciated. Further evidence comes from some foreign visitors, who wrote: "I wish I had paid more attention to math class!" indicating the fact that the exhibition has managed to somehow overcome the misconception regarding mathematics and bring in the foreground some of the fundamental aspects of the teaching of mathematics (but not only), such as curiosity, wonder and the transmission of important ideas and not only tedious techniques. In this regard, a visitor wrote that he/she discovered a new world. In this sense, the initial disorientation was overcome in an unexpected and positive way because the exhibition succeeded in unveiling the potential of mathematics, which previously seemed unattainable. Another witness is grateful for the many ideas and so many hints that would be used at school for events on the subject of mathematics. There have also been some criticism, sometimes because of the difficulty in understanding certain topics. For this, as we discussed before, the guides were really fundamental because they played a role as a cultural mediator. The comments regarding some guides have been very positive, as, for example, "The exhibition is very beautiful but the enthusiasm of our guide was exceptional". Finally, we would like to quote a comment by a visitor, probably very young, that after the visit wrote that the exhibition was better than \textit {Gardaland} (an amusement park)!

\section{Conclusions}\label{final}

From what has been written so far, the exhibition \textit{MadeInMath} might seem so multifaceted that the various threads fail to help us come out of the maze in which mathematics could make us lose. On the other hand, all these threads are unraveled from one skein containing thousands of answers to the question that marked the beginning of the design of the exhibition: what is the profession of the mathematician? A question that has stimulated visitors to find answers as well, because there is certainly not only one!

\end{document}